\documentclass[pdflatex,sn-mathphys-num]{sn-jnl}

\usepackage{graphicx}%
\usepackage{multirow}%
\usepackage{amsmath,amssymb,amsfonts}%
\usepackage{amsthm}%
\usepackage{mathrsfs}%
\usepackage[title]{appendix}%
\usepackage{xcolor}%
\usepackage{textcomp}%
\usepackage{manyfoot}%
\usepackage{booktabs}%
\usepackage{algorithm}%
\usepackage{algorithmicx}%
\usepackage{algpseudocode}%
\usepackage{listings}%

\usepackage{subcaption}
\usepackage{tikz}
\usepackage{siunitx}
\usepackage{hyperref}
\usepackage{bbm}

\DeclareMathOperator{\tr}{tr}
\newcommand{\R}{\mathbb{R}}
\newcommand{\N}{\mathbb{N}}

\theoremstyle{thmstyleone}%
\newtheorem{theorem}{Theorem}[section]
\theoremstyle{thmstyletwo}%

\theoremstyle{thmstylethree}%
\newtheorem{definition}[theorem]{Definition}%

\begin{document}

\title[Multi-Criteria Inverse Robustness in Radiotherapy Planning Using Semidefinite Programming]{Multi-Criteria Inverse Robustness in Radiotherapy Planning Using Semidefinite Programming}


\author[1]{\fnm{Jan} \sur{Schr\"oeder}}\email{jan.schroeder@itwm.fraunhofer.de}

\author[2]{\fnm{Yair} \sur{Censor}}\email{yair@math.haifa.ac.il}
\equalcont{These authors contributed equally to this work.}

\author*[1]{\fnm{Philipp} \sur{S\"uss}}\email{philipp.suess@itwm.fraunhofer.de}
\equalcont{These authors contributed equally to this work.}

\author[1]{\fnm{Karl-Heinz} \sur{K\"ufer}}\email{karl-heinz.kuefer@itwm.fraunhofer.de}
\equalcont{These authors contributed equally to this work.}

\affil*[1]{\orgdiv{Optimization Department}, \orgname{Fraunhofer ITWM}, \orgaddress{\street{Fraunhofer-Platz 1}, \city{Kaiserslautern}, \postcode{67663}, \country{Germany}}}

\affil[2]{\orgdiv{Department of Mathematics}, \orgname{University of Haifa}, \orgaddress{\street{Abba Khoushy Ave 199}, \city{Haifa}, \postcode{3498838}, \country{Israel}}}


\abstract{Radiotherapy planning naturally leads to a multi-criteria optimization problem which is subject to different sources of uncertainty. In order to find the desired treatment plan, a decision maker must balance these objectives as well as the level of robustness towards uncertainty against each other. This paper showcases a quantitative approach to do so, which combines the theoretical model with the ability to deal with practical challenges. To this end, the uncertainty, which can be expressed via the so-called dose-influence matrix, is modelled using interval matrices. We use inverse robustness to introduce an additional objective, which aims to maximize the volume of the uncertainty set. A multi-criteria approach allows to handle the uncertainty while keeping appropriate values of the other objective functions. We solve the resulting quadratically constrained quadratic optimization problem (QCQP) by first relaxing it to a convex semidefinite problem (SDP) and then reconstructing optimal solutions of the QCQP from solutions of the SDP.}




\maketitle

\section{Introduction}
In the fight against cancer, radiotherapy plays a crucial role. About half of all cancer patients receive radiation as a part of their treatment \cite{cancerorg}. Planning radiotherapy treatments involves multiple conflicting objectives; in particular, dose maximization in the tumor and dose minimization in organs at risk (OARs) are typically competing forces to balance. Numerical optimization of most radiotherapy plans is typically done using a linear model for the dose deposited in the body as a function of the radiation. The total dose $d_i$ absorbed in voxel $i$ is calculated by  $d=Dx$, where $D$ is the dose-influence matrix, $d$ is the dose vector and $x$ is the vector of fluences (radiation). 

One complication to treatment planning are uncertainties \cite{unkelbach}. Fortunately, many uncertainties (for example organ motion, setup errors, range errors) can be represented in the matrix $D$ \cite{unkelbach}. In this paper we consider interval uncertainty for the entries $D_{ij}$ of $D$ and enhance it with concepts from inverse robustness \cite{Berthold2024} . This allows to adjust the level of robustness (e.~g. ``$50\%$ of scenarios should be covered'') adaptively to the desired amount while keeping the other objectives Pareto-optimal. In this way a decision maker can navigate the appropriate trade-off between tumor dose, OAR dose and robustness. 

We develop and present a method that, using the above techniques, allows to efficiently compute the Pareto-optimal decisions for this problem. To numerically handle the resulting large-sized problem, we apply a voxel clustering method first suggested in \cite{scherrer}, which can significantly reduce the problem size and computational effort, while yielding qualitatively approximately equivalent solutions.

The paper is laid out as follows: We begin with an overview of established methods and the tools that are used in the sequel in Section \ref{prel}. Starting from the problem formulation and the introduction of scalarizations and uncertainties via interval matrices, we move on to the SDP-relaxation and end with a description of a clustering method to reduce the problem size. These steps lead to a relaxed semidefinite optimization problem. We introduce a reconstruction method, which allows to obtain solutions for the non-relaxed problem in Section \ref{recmeth}. Finally, we show how this method can be interwoven with Pareto-front navigation for better decision making in Section \ref{paretonavmeth}. We apply our method to a real-world case and compare the outcome with results from established methods for uncertainty handling in Section \ref{compmeth}. We conclude the paper in Section \ref{conclusion}, where we recapitulate the most important findings.

\subsection{Context and Previous Work}
In this work we make use of a linear problem formulation like the one proposed in TROTS \cite{trots}. TROTS is a test set of clinical cases for both photon and proton therapy planning. It includes data for multiple cases like liver, prostate or head-and-neck. It features a voxel-wise formulation, which can be adapted to incorporate constraints on (weighted) means or average constraint violation while maintaining the linear character \cite{gorissen}. 

Earlier works employed non-linear models (\cite{bokrantz,unkelbach}), but we make use of the structural properties that come with a linear formulation and make a similar argument for the introduction of uncertainty. By using interval arithmetic (\cite{rohn}) to model and handle uncertainty, we keep the model simple and obtain a quadratically constrained quadratic problem (QCQP). An approach of using interval matrices and robustness levels has recently been proposed in \cite{moreno}, but there the authors make multiple computations for different levels of robustness. In our work, we enhance this approach by using inverse robustness (\cite{Berthold2024}), which allows for a better trade-off analysis and Pareto-front exploration. 

Other approaches model uncertainty in terms of scenarios (\cite{bokrantzfredriksson}) or via covariances (\cite{wahl}) and make use of stochastic methods (\cite{bangert}). Our results on the formulation of the QCQP and its SDP relaxation \cite{vandenberghe} fall in the category of similar results on the recovery of solutions from relaxations as in \cite{kim,kimizuka} or \cite{azuma}.

\section{Preparing the Ground}\label{prel}
In this section we set up the relevant problems and tools whose accumulated adoption constitutes our new navigation method. We begin with the radiotherapy formulation as a linear problem (LP), induced by the multi-criteria approach, which after the introduction of uncertainty and interval matrices becomes a quadratically constrained quadratic problem (QCQP). We relax this QCQP to obtain a semidefinite problem (SDP), which we cluster, to reduce its size, and to which we can then apply the reconstruction method.

\subsection{The Radiotherapy Problem}\label{imrtProblem}
Denote by $\tilde{m},\tilde{n}\in\N$ the number of voxels and the number of beamlets, respectively. Denote by $D\in\R^{\tilde{m}\times \tilde{n}}$ the dose-influence matrix. Let $\mathcal{S}$ be the set of all structures that have been identified within the irradiated part of the body. For each $s\in \mathcal{S}$ the set $I(s)\subseteq \{1,2,\dots,\tilde{m}\}$ is a list of voxels that are contained in structure $s$. Note that a voxel can contain parts of more than one structure and so the sets $I(s)$ need not have an empty intersection.

Further, we choose subsets $\mathcal{O}$ and $\mathcal{C}$ of $\mathcal{S}$, on which the dose is to be optimized and on which the dose has constraints, respectively. Again, the intersection $\mathcal{O}\cap \mathcal{C}$ need not be empty, but it can be. Further, for each $s\in \mathcal{C}$, let $b_s^\text{lb},b_s^\text{ub}\in\R$ be the lower and upper dose bounds, respectively. The full multi-criteria problem can then be stated as
\begin{subequations}\label{eq:imrtprob}
	\begin{align}
		\min_{x\in\R_\geq^{\tilde{n}},z\in\R^{\vert \mathcal{O}\vert}}& \begin{pmatrix}
			z_1\\\vdots\\z_{\vert \mathcal{O}\vert}
		\end{pmatrix}\\
		\text{s.t. }&b^\text{lb}_s\leq\sum_{j=1}^{\tilde{n}} D_{ij}x_j\leq b^\text{ub}_s\quad\forall i\in I(s),s\in \mathcal{C}\\\label{eq:imrtobjepi}
		&\sum_{j=1}^{\tilde{n}}D_{ij}x_j\leq z_s\quad\forall i\in I(s),s\in \mathcal{O},
	\end{align}
\end{subequations}

where the minimization is to be understood in the Pareto-sense (which we define next) and $\vert \mathcal{O}\vert$ stands for the cardinality of $\mathcal{O}$. This problem formulation agrees with the one proposed in TROTS \cite{trots}.

\subsection{Scalarizations in Multi-Criteria Optimization}\label{scalarizations}
Let $\mathcal{X}\subseteq\R^n$ and let $f\colon \mathcal{X}\to\R^k$. A general multi-criteria optimization problem is
\begin{align}
	\min_{x\in \mathcal{X}}\ f(x).\label{eq:mco}
\end{align}

Such a multi-criteria problem requires an adjusted notion of optimality as, in general, not all components of the objective function will attain their optimum for the same value of $x$. We use Pareto-optimality for this problem, a notion which dates decades ago and we briefly review below:

\begin{definition}[{\cite[Definition 2.1, 2.24]{ehrgott}}]
	\begin{enumerate}[(i)]
		\item A feasible solution $x^*\in \mathcal{X}$ of \eqref{eq:mco} is called Pareto-optimal, a Pareto-point, or efficient, if there is no $x\in \mathcal{X}$ such that $f(x)$ dominates $f(x^*)$, that is $f(x)\leq f(x^*)$ and $f(x)\neq f(x^*)$. In this case, $f(x^*)$ is called non-dominated or a Pareto-point. The set of all efficient solutions $x^*\in \mathcal{X}$ is denoted by $\mathcal{X}_E$. The set of all non-dominated points $y^*=f(x^*)\in \mathcal{Y}$ is denoted by $\mathcal{Y}_N$. We also call $\mathcal{Y}_N$ the Pareto-front.
		\item A feasible solution $x^*\in \mathcal{X}$ of \eqref{eq:mco} is called weakly Pareto-optimal, a weak Pareto-point or weakly efficient, if there is no $x\in \mathcal{X}$ such that $f(x)<f(x^*)$. In this case, $f(x^*)$ is called weakly non-dominated or a weak Pareto-point. The set of all weakly efficient solutions $x^*\in \mathcal{X}$ is denoted by $\mathcal{X}_{wE}$. The set of all weakly non-dominated points $y^*=f(x^*)\in \mathcal{Y}$ is denoted by $\mathcal{Y}_{wN}$.
		\item A feasible solution $x^*\in \mathcal{X}$ of \eqref{eq:mco} is called strictly Pareto-optimal, a strict Pareto-point or strictly efficient, if there is no other $x\in \mathcal{X}$ such that $f(x)\leq f(x^*)$. In this case, $f(x^*)$ is called strictly non-dominated or a strict Pareto-point. The set of all strictly efficient solutions $x^*\in \mathcal{X}$ is denoted by $\mathcal{X}_{sE}$. The set of all strictly non-dominated points $y^*=f(x^*)\in \mathcal{Y}$ is denoted by $\mathcal{Y}_{sN}$.
	\end{enumerate}		
\end{definition}

Another useful notion is $\epsilon$-efficiency:

\begin{definition}[{\cite[Definition 2.1]{liu}}]
	Let $\epsilon\geq0$. A feasible solution $x^*\in\mathcal{X}$ is called $\epsilon$-Pareto-optimal or $\epsilon$-efficient, if there is no $x\in\mathcal{X}$ which fulfills $f(x)\leq f(x^*)-\epsilon\mathbbm{1}$ and $f(x)\neq f(x^*)-\epsilon\mathbbm{1}$, where $\mathbbm{1}$ is the vector of all ones of appropriate dimension. Similarly, $x^*$ is called weakly $\epsilon$-Pareto-optimal or weakly $\epsilon$-efficient, if there is no $x\in\mathcal{X}$ with $f(x)< f(x^*)-\epsilon\mathbbm{1}$.
\end{definition}
Note that for $\epsilon=0$ this definition coincides with the classical definition of Pareto-optimality.

In order to find Pareto-optimal points we use scalarizations, i.~e., methods that reduce the multi-criteria problem to a single-criteria problem that finds efficient points of \eqref{eq:mco}. The first method, which is particularly useful in a convex setting, is the weighted sum scalarization. Here, each component objective $f_i$ is assigned a weight $w_i$ and the weighted sum is optimized. Usually the weights $w_i$ are chosen to be non-negative or positive in order to ensure efficiency of the resulting point $x^*$ \cite[Theorem 4.1]{ehrgott}. The weighted sum optimization problem with weight $w\in\R^k$ associated with the multi-criteria problem \eqref{eq:mco} then is
\begin{align}
	\min_{x\in \mathcal{X}}&\ w^Tf(x).\label{eq:wsopt}
\end{align}

Another scalarization technique, the $\epsilon$-constraint method, introduces a vector parameter $\epsilon\in\R^k$, which is used to impose bounds on all but one of the objective functions, say the $j$-th ($j\in\{1,2,\dots,k\}$). This leads to the following single-criteria problem:
\begin{subequations}\label{eq:epsconstr}
	\begin{alignat}{3}
		\min_{x\in \mathcal{X}} &\ f_j(x) &\\
		\text{s.t.} &\ f_i(x)\leq\epsilon_i,\quad & i=1,2,\dots,k,\ i\neq j.
	\end{alignat}
\end{subequations}
Again, we can conclude efficiency for \eqref{eq:mco} from optimality for the scalarization \eqref{eq:epsconstr} (see, e.~g., \cite[Proposition 4.3]{ehrgott}). Note, that the $\epsilon$-constraint method also works in a non-convex setting.

We will use these techniques to scalarize our multi-criteria problem. We will see that applying the weighted sum method to the (relaxed) convex problem leads to much fewer Pareto-points that need to be computed.

\subsection{Uncertainty and Inverse Robustness}\label{uncir}

For ease of notation we will consider the problem 

\begin{subequations}\label{eq:mylp}
	\begin{align}
		\min_{x\in\R^n}\ &Gx\\
		\text{s. t. } & Ax\leq b\\
		&\ell\leq x\leq u,
	\end{align}
\end{subequations}

where $k,m,n\in\N$, $A\in\R^{m\times n}$, $b\in\R^m$ and $G\in\R^{k\times n}$. Clearly, problem \eqref{eq:imrtprob} can be stated in this notation.

Uncertainty plays an important role in radiotherapy. All uncertainty within the patients anatomy can be incorporated into the dose-influence matrix $D$ \cite{unkelbach}. For our purposes, this means considering uncertainty in the matrix $A$. While $G$ may also contain uncertainty, we can reformulate uncertainties in the objectives to be contained in the constraints instead, by using the well-known epigraph reformulation \cite[problem (4.11)]{boyd}. 

The formulation \eqref{eq:mylp} can accomodate a problem statement like in \eqref{eq:imrtprob}, where the objectives are constrained by the auxiliary variables $z_i$. The bounds $\ell$, $u$ and the right-hand side $b$ are prescribed by the physician, but uncertain $\ell$, $u$ and $b$ can also be treated within the following framework. 

Let $U$ be the uncertainty set, i.~e., the set of all $W:=(A,b)$. The classical worst-case formulation \cite{bental} of an ``uncertain inequality'' $Ax\leq b$ with data $A,b$ is
\begin{align}\label{eq:bentalwc}
	Ax\leq b\quad\forall (A,b)\in U.
\end{align}

Often times, handling uncertainty in this way leads to results that are too conservative. This observation gives rise to the concept of ``inverse robustness''. Here, an additional objective is to maximize a utility function $\vartheta\colon\mathcal{W}\to\R$, which is defined over a collection of subsets of $U$, i.~e., $\mathcal{W}\subseteq 2^U$. Namely, we add

\begin{align}
	\max_{W\in\mathcal{W}} \vartheta (W)
\end{align}

to the objectives and modify the classical worst-case constraint \eqref{eq:bentalwc} to read

\begin{align}\label{eq:invrobconstr}
	Ax\leq b\quad\forall(A,b)\in W.
\end{align}

Maximizing a utility function over the uncertainty set in this way allows to explore trade-offs between the magnitude of deviation from a nominal scenario and other objective functions. Below, we embed this approach into the framework of interval matrices. This leads to a formulation that has the relevant properties to ensure the existence of optimal solutions (\cite{Berthold2024}) and the structure to efficiently calculate them.

\subsection{Interval Matrices}\label{intmat}

In this subsection we use interval matrices to model uncertainty in the radiotherapy planning problem and subsequently apply the concept of inverse robustness to obtain a quadratically constrained quadratic problem (QCQP).

\begin{definition}[{\cite[Section 2.5]{rohn}}]
	Let $\underline{A},\overline{A}\in\mathbb{R}^{m\times n}$ be real matrices. An interval matrix is the set of matrices
	\begin{equation}
		\boldsymbol{A}=[\underline{A},\overline{A}]=\{A\in\mathbb{R}^{m\times n}\vert\ \underline{A}\leq A\leq\overline{A}\},
	\end{equation}
	where the inequalities are to be understood component-wise. Observe that we use boldface symbols for inter val quantities. The center matrix $A_c=\frac{1}{2}(\overline{A}+\underline{A})$ and the offset matrix $A_\delta=\frac{1}{2}(\overline{A}-\underline{A})$ are also used, as needed for the given context.
	
	An interval vector is an interval matrix with only one column: Let $\underline{b},\overline{b}\in\mathbb{R}^m$. An interval vector is the set
	\begin{equation}
		\boldsymbol{b}=[\underline{b},\overline{b}]=\{b\in\mathbb{R}^{m}\vert\ \underline{b}\leq b\leq\overline{b}\}.
	\end{equation}
\end{definition}

In the following let $\boldsymbol{A}$ be an $m\times n$ interval matrix and $\boldsymbol{b}$ an $m$-dimensional interval vector. We consider the system 
\begin{equation}
	\boldsymbol{A}x\leq \boldsymbol{b}
\end{equation}
as the family of all systems
\begin{equation}
	Ax\leq b\text{, with } A\in\boldsymbol{A},b\in\boldsymbol{b}.
\end{equation}

A solution is defined as follows:

\begin{definition}[{\cite[Section 2.13]{rohn}}]
	\begin{enumerate}[(i)]
		\item The system $\boldsymbol{A}x\leq \boldsymbol{b}$ is called strongly solvable if for every $A\in\boldsymbol{A}, b\in\boldsymbol {b}$ there is some $x\in\mathbb{R}^n$ with $Ax\leq b$. 
		
		\item The vector $x$ is called a strong solution of $\boldsymbol{A}x\leq\boldsymbol{b}$ if $Ax\leq b$ for all $A\in\boldsymbol{A}, b\in\boldsymbol {b}$.
	\end{enumerate}
\end{definition}

Note that strong solvability only needs some individual solution $x$ depending on the actual realizations $A$ and $b$, whereas a strong solution will solve the interval inequality regardless of the realization. Nevertheless the following theorem holds:

\begin{theorem}[{\cite[Theorem 2.24]{rohn}}]
	If $\boldsymbol{A}x\leq\boldsymbol{b}$ is strongly solvable then it has a strong solution.
\end{theorem}

For the proof we refer to \cite{rohn}. The next theorem gives a verifiable condition for a vector $x$ to be a strong solution of $\boldsymbol{A}x\leq\boldsymbol{b}$. Here, $x_+=\max(x,0)$ $x_-=\max(-x,0)$ and $\vert x\vert$ denote the positive part, the negative part and the absolute value of $x$, respectively, and all operators are to be understood component-wise:

\begin{theorem}[{\cite[Theorem 2.25]{rohn}}]\label{IntMatLinear}
	Let $\boldsymbol{A}\subseteq\R^{m\times n}$ be an interval matrix, $\boldsymbol{b}\subseteq\R^m$ be an interval vector and $x\in\R^n$. The following assertions are equivalent:
	\begin{enumerate}
		\item $x$ is a strong solution of $\boldsymbol{A}x\leq\boldsymbol{b}$
		\item $x$ satisfies
		\begin{align}\label{eq:rohn2}
			A_cx-b_c\leq -A_\delta\vert x\vert-b_\delta
		\end{align}
		
		\item $x=x_+-x_-$, where $x_+, x_-\in\R^n$ satisfy 
		\begin{subequations}\label{eq:rohn3}
			\begin{align}
				\overline{A}x_+-\underline{A}x_-&\leq\underline{b},\label{eq:rohn3a}\\
				x_+, x_-&\geq 0.
			\end{align}
		\end{subequations}
	\end{enumerate}
\end{theorem}

To embed interval matrices within the context of inverse robustness, we set up $\mathcal{W}$ and $\vartheta$ as follows:	

The uncertainty set $U$ is the set of all pairs $(A,b)$, i.e., 

\begin{align}\label{eq:invrobsetupbegin}
	U:=\{(A,b)\colon A\in\boldsymbol{A},b\in\boldsymbol{b}\}=[\underline{A},\overline{A}]\times[\underline{b},\overline{b}].
\end{align}

We introduce a real parameter $r$, which will be used as a scaling variable for the inverse robustness, and consider sets of the form

\begin{align}
	W(r):=[A_{c}-rA_{\delta},A_{c}+rA_{\delta}]\times[b_{c}-rb_{\delta},b_{c}+rb_{\delta}].
\end{align}
We define  
\begin{align}
	\mathcal{W}:=\{W(r)\colon r\in[0,1]\}
\end{align} 
and the utility function as
\begin{align}\label{eq:invrobsetupend}
	\vartheta(W(r)):=r.
\end{align}

By defining the setup and parameters as in \eqref{eq:invrobsetupbegin}-\eqref{eq:invrobsetupend}, we are perfectly within the framework from \cite{Berthold2024} and can use the theorems on the existence of optimal solutions there.

In this case we apply  Theorem \ref{IntMatLinear} to the semi-infinite constraint \eqref{eq:invrobconstr} and obtain from \eqref{eq:rohn3a}

\begin{align}
	(A_{c}+rA_{\delta})x_{+}-(A_{c}-rA_{\delta})x_{-}\leq b_{c}-rb_{\delta}.
\end{align}	

Since $r$ is an optimization parameter, this is a quadratic constraint. After these steps, the full QCQP is

\begin{subequations}\label{eq:myqcqp}
	\begin{align}
		\min_{x_+,x_-\in\R_\geq^n,\ r\in[0,1]}&\begin{pmatrix} G(x_+-x_-)\\ -r\end{pmatrix}\\
		\text{s.t. } & A_cx_+ -A_cx_- +rA_\delta x_+ +rA_\delta x_-+rb_\delta -b_c\leq 0\label{eq:myqcqpcon}\\
		&\ell\leq x_+-x_-\leq u.\label{eq:myqcqpcon2}
	\end{align}
\end{subequations}

For given fixed $r\geq0$ we set 
\begin{align}
	\overline{A}_r&=A_c+rA_\delta\\
	\underline{A}&=A_c-rA_\delta
\end{align}
and analogously for $\overline{b}_r$ and $\underline{b}_r$. We call the problem
\begin{subequations}
	\begin{align}
		\min_{x_+,x_-\in\R_\geq^{n}}&G(x_+-x_-)\\
		\text{s.t. } & \overline{A}_rx_+-\underline{A}_rx_-\leq\underline{b}_r,
	\end{align}
\end{subequations}
the (multi-criteria) LP of robustness level $r$ associated with the QCQP \eqref{eq:myqcqp}.

\subsection{SDP-Relaxation}\label{sdprel}
In general, QCQPs are non-convex and NP-hard to solve. Hence, we make use of a relaxation technique, the SDP-relaxation, to obtain a convex optimization problem, which can be solved efficiently. We first describe the SDP-relaxation for general QCQPs of the form
\begin{subequations}\label{eq:qcqp}
	\begin{align}
		\min_{x\in\R^N}&\ x^TQ_0x+2q_0^Tx+\gamma_0\\
		\text{s.t.} &\  x^TQ_ix+2q_i^Tx+\gamma_i\leq 0,\ i=1,2,\dots,m
	\end{align}
\end{subequations}
with $m,N\in\N$, square matrices $Q_i\in\R^{N\times N}$, vectors $q_i\in\R^N$ and scalars $\gamma_i\in\R$ ($i=0,\dots,m$). Without loss of generality we assume that all $Q_i$ are symmetric, as we can otherwise replace them with the symmetric matrix $\frac{1}{2}(Q_i+Q_i^T)$ without changing the overall form of the objective function or constraint.

Define for $i\in\{0,\dots,m\}$ the matrices
\begin{align}
	M_i=\begin{pmatrix}
		\gamma_i&q_i^T\\
		q_i&Q_i
	\end{pmatrix}
\end{align}
and 
\begin{align}
	M_{m+1}=\begin{pmatrix}1&0_N^T\\0_N&0_{N\times N}\end{pmatrix}.
\end{align}
Then problem \eqref{eq:qcqp} can be stated as 
\begin{subequations}\label{eq:nqcqp}
	\begin{align}
		\min_{z\in\R^{N+1}}&\ z^TM_0z\\
		\text{s.t.} &\  z^TM_iz\leq 0,\ i=1,2,\dots,m\\
		&z^TM_{m+1}z=1,
	\end{align}
\end{subequations}
where $z=(z_0,x)$ contains $x$ and an additional variable $z_0$. Clearly, due to the additional constraint from $M_{m+1}$ we get $z_0=\pm1$ and we see that $z=(z_0,x)$ is feasible for \eqref{eq:nqcqp} if and only if $z_0x$ is feasible for \eqref{eq:qcqp}. Further, we find that $(z_0,x)$ and $(-z_0,-x)$ give the same value for the objective function and for each constraint, so that we can, without loss of generality, assume that $z_0=1$.

We rewrite problem \eqref{eq:nqcqp} once again to get
\begin{subequations}
	\begin{align}
		\min_{z\in\R^{N+1}}&\ \langle M_0,zz^T\rangle\\
		\text{s.t.} &\  \langle M_i,zz^T\rangle\leq 0,\ i=1,2,\dots,m\\
		&\langle M_{m+1},zz^T\rangle=1,
	\end{align}
\end{subequations}
where $\langle\cdot,\cdot\rangle$ denotes the matrix inner product, i.e., $\langle A,B\rangle=\tr(A^TB)$ where $\tr$ stands for the matrix trace operator. It is easy to see that the matrix $zz^T$ is symmetric, positive semidefinite and has rank $1$ since it is generated by a singular vector. We introduce a matrix $Z\in\R^{(N+1)\times (N+1)}$ that replaces $zz^T$ and drop the so-called rank-1-constraint that $Z$ be of rank $1$. In this way, we obtain a relaxation of the original problem, the well-known SDP-relaxation:
\begin{subequations}\label{eq:sdprel}
	\begin{align}
		\min_{Z\in S_+^{N+1}}&\ \langle M_0,Z\rangle\\
		\text{s.t.} &\  \langle M_i,Z\rangle\leq 0,\ i=1,2,\dots,m\\
		&\langle M_{m+1},Z\rangle=1.
	\end{align}
\end{subequations}
Here, $S_+^{N+1}$ is the cone of symmetric and positive semidefinite matrices in $\R^{(N+1)\times (N+1)}$. This relaxation is convex and the only non-convexity in the QCQP stems from the rank-1-constraint $Z=zz^T$. The constraint $Z\in S_+^{N+1}$ is convex conic in nature and all other constraints and the objective function are linear in $Z$. The following theorem relates optimality in the SDP to optimality in the QCQP:

\begin{theorem}[{\cite[Section 2]{vandenberghe}}]\label{thmrank1opt}
	Let $Z$ be optimal for \eqref{eq:sdprel}. If $Z$ is of rank 1, i.~e., there is some $z$ with $Z=zz^T$, then $z$ is optimal for \eqref{eq:nqcqp}.
\end{theorem}

Note that optimality of $z=(z_0,x)$ in Theorem \ref{thmrank1opt} further implies optimality of $x$ in \eqref{eq:qcqp} as these problems can equivalently be reformulated into each other.

Let us return to our radiotherapy treatment planning problem \eqref{eq:myqcqp}. In order to make use of the SDP-relaxation, we phrase it in terms of matrices $M_i$. To this end, define 

\begin{align}\label{eq:sdpconmat}
	M_i=\frac{1}{2}
	\begin{pmatrix}
		-2b_c^i & (A_c^i)^T & -(A_c^i)^T & b_\delta^i\\
		A_c^i & 0_{n\times n} & 0_{n\times n} & A_\delta^i\\
		-A_c^i & 0_{n\times n} & 0_{n\times n} & A_\delta^i\\
		b_\delta^i & (A_\delta^i)^T & (A_\delta^i)^T & 0
	\end{pmatrix},
\end{align}
for $i=1,2,\dots,m$, where $(A_c^i)^T$ denotes the $i$-th row of $A_c$ and analogously for $(A_\delta^i)^T$ and $(G^i)^T$. Similarly, for the $k+1$ objective functions set

\begin{align}\label{eq:sdpobjmat}
	M_0^i=\frac{1}{2}
	\begin{pmatrix}
		0 & (G^i)^T & -(G^i)^T & 0\\
		G^i & 0_{n\times n} & 0_{n\times n} & 0_n\\
		-G^i & 0_{n\times n} & 0_{n\times n} & 0_n\\
		0 & 0_n^T & 0_n^T & 0
	\end{pmatrix},
\end{align}
for $i\in\{1,2,\dots,k\}$ as well as 
\begin{align}\label{eq:sdprmat}
	M_0^{k+1}=\frac{1}{2}
	\begin{pmatrix}
		0 & 0_n^T & 0_n^T & -1\\
		0_n & 0_{n\times n} & 0_{n\times n} & 0_n\\
		0_n & 0_{n\times n} & 0_{n\times n} & 0_n\\
		-1 & 0_n^T & 0_n^T & 0
	\end{pmatrix}.
\end{align}
The matrix $M_{m+1}$ remains the same. Note that the constraints $\ell\leq x_+-x_-\leq u$ are not subject to uncertainty and can be modelled as $\ell\leq Ix_+-Ix_-\leq u$, where $I$ is the identity matrix. In this way, we can handle these bounds in the same way as the other constraints and, therefore, we will leave them out in the following analysis. Bear in mind though that they remain part of the optimization problem and we require them here and there. 

Note that we have some additional information that we can use. Namely, with $Z_{11}=1$ being a ``constant variable'', we can interpret the rest of the first row (or column) to represent the linear terms $\hat{x}^T=(x_+,x_-,r)^T$. Similarly, the remaining block matrix represents all the quadratic variables, i.e., the products $\hat{x}_i\hat{x}_j$, which have been condensed into a single matrix variable, which we call $X$, whose $ij$-th element is $X_{ij}=\hat{x}_i\hat{x}_j$. We get 
\begin{align}
	Z\approx\begin{pmatrix} 1 & \begin{pmatrix} x_+&x_-&r\end{pmatrix}\\\begin{pmatrix}x_+\\x_-\\r\end{pmatrix}& X\end{pmatrix}=\begin{pmatrix}1&\hat{x}^T\\\hat{x}& X\end{pmatrix}.
\end{align}
This intuition guides us throughout this section and we will see that, in this particular case, it is in fact quite helpful. For example, it allows us to impose the constraint $Z\geq0$, since we know that all variables $x_+,x_-,r$ are non-negative and so their products must be too. 

From $\ell\leq x\leq u$ we also see (with a slight abuse of notation) that $\hat{x}_i^2\approx X_{ii}\leq \max\{\tilde{u}_i^2,\tilde{\ell}_i^2\}$ should hold for any meaningful solution. Here, we extended the bounds $\ell$ and $u$ to the appropriate dimensions of variables: $\tilde{\ell}=(\ell^T,\ell^T,0)^T$ and $\tilde{u}=(u^T,u^T,1)^T$, respectively. 

In other words, introducing all these constraints into the relaxed SDP problem helps us reduce the feasible set, removing any solutions that are not meaningful anyway, without ruling out solutions that are feasible for the original QCQP, since the added constraints are implicitly contained in there already.

The resulting multi-criteria semidefinite optimization problem can be scalarized using scalarization techniques like the weighted-sum method, desribed in Section \ref{scalarizations} above. For example with weights $w_i\geq0$ we can set $M_0=\sum_{i=1}^{k+1}w_iM_0^i$. Due to the convexity of the SDP, we can find all Pareto-optimal points like this and approximate the Pareto-front to the desired accuracy. 

Note that a solution of the SDP \eqref{eq:sdprel} will, in general, be infeasible for the QCQP (i.~e., $Z$ will not be of rank $1$ and no $z$ with $Z=zz^T$ exists) and it is not obvious how to construct a ``nearby'' rank-$1$-solution from an SDP-solution $Z$. Nevertheless, due to the particular structure of our SDP, we can formulate such a reconstruction method for our particular problem class and we do so in Section \ref{recmeth}.

\subsection{Clustering}\label{clustering}
Before getting into the details of the reconstruction method, we address the concern about dimensionality. Radiotherapy problems tend to be large. Our above reformulations make the problem even larger, for example $n$ variables in $x$ become a $(2n+2)\times (2n+2)$-matrix. Further, modern primal-dual solvers for SDPs like MOSEK \cite{mosek} construct the dual problem too, which then contains additional variables for each constraint. Storing these variables as standard 16-bit floats then quickly exceeds the memory on today's off-the-shelf computers. Consequently, we need a method to reduce the problem size, or, more precisely, to reduce the number of constraints on the individual voxels.

To do this, we cluster similar voxels together to create multiple super-voxels. We perform two steps for each individual structure (so that two different structures can not be clustered together into the same super-voxel). In the first step, we perform a statistical analysis of each row of the dose-influence matrix. We calculate the statistical mean and the statistical variance within each row and accumulate them in a list of size $2\times m$. 

In the second step we apply the $k$-means algorithm \cite{arthur} to this list to obtain a number of $K$ clusters of voxels that show a similar behaviour. We combine these voxels into a single super-voxel by taking the mean over the columns in the dose-influence matrix (only considering entries from the individual cluster). In this way, we get $K$ clusters, which we combine into $K$ super-voxels, shrinking the dose-influence matrix from size $m\times n$ to $K\times n$. 

In \cite{scherrer} it was shown that using this method significantly reduces the problem size and increases computation speed while effectively leaving the resulting treatment plans at a comparable quality. This gives us confidence that this method will allow us to solve very large problems to an acceptable quality. A similar method exists for clustering of variables \cite{suess}.

\section{Reconstruction Method}\label{recmeth}
In order to reconstruct a feasible solution of the QCQP, assume that we found a solution 
\begin{align}
	Z^\text{SDP}=\begin{pmatrix}1& \begin{pmatrix}x_+^\text{SDP}&x_-^\text{SDP}&r^\text{SDP}\end{pmatrix}\\\begin{pmatrix}x_+^\text{SDP}\\x_-^\text{SDP}\\r^\text{SDP}\end{pmatrix}& X^\text{SDP}\end{pmatrix}
\end{align}
of the SDP. We set $x_+^\text{QCQP}=x_+^\text{SDP}$ and $x_-^\text{QCQP}=x_-^\text{SDP}$ and choose $r^\text{QCQP}$ as the solution of 

\begin{subequations}\label{eq:maxr}
	\begin{align}
		\min_{r\in[0,1]} &-r\\
		\text{s.t. }&(A_\delta x_+^\text{SDP} +A_\delta x_-^\text{SDP}+b_\delta)r + A_cx_+^\text{SDP} -A_cx_-^\text{SDP}  -b_c\leq0.\label{eq:maxrcon}
	\end{align}
\end{subequations}

This is a linear problem with a single optimization variable that can be solved within milliseconds, even for large $m,n$. It is easy to see that $r=0$ is always feasible for \eqref{eq:maxr} due to the feasibility of $Z^\text{SDP}$ for \eqref{eq:sdprel}. Hence, the feasible set of \eqref{eq:maxr} is non-empty and compact and an optimal solution always exists. A number of theorems can be proven for $(x_+^\text{QCQP}, x_-^\text{QCQP}, r^\text{QCQP})$, which we only state here without proof.

\begin{theorem}[Efficiency of the projected point, {\cite[Theorem 5.8]{schroeder}}]\label{thmoptstr}
	Let $(x_+^\text{QCQP},x_-^\text{QCQP},r^\text{QCQP})^T$ be obtained from an SDP-point
	\begin{align}\label{eq:sdppoint}
		Z^\text{SDP}=\begin{pmatrix}1& \begin{pmatrix}x_+^\text{SDP}&x_-^\text{SDP}&r^\text{SDP}\end{pmatrix}\\\begin{pmatrix}x_+^\text{SDP}\\x_-^\text{SDP}\\r^\text{SDP}\end{pmatrix}& X^\text{SDP}\end{pmatrix}.
	\end{align} If $(x_+^\text{QCQP},x_-^\text{QCQP})$ is a strictly efficient solution for the LP of robustness level $r^\text{QCQP}$,
	then $(x_+^\text{QCQP},x_-^\text{QCQP},r^\text{QCQP})^T$ is efficient for the QCQP \eqref{eq:myqcqp}.
\end{theorem}

\begin{theorem}[Weak $\epsilon$-efficiency of the projected point, {\cite[Theorem 5.9]{schroeder}}]\label{thmoptw}
	Let $\epsilon\geq0$ and let the vector $(x_+^\text{QCQP},x_-^\text{QCQP},r^\text{QCQP})^T$ be obtained from an SDP-point as in\eqref{eq:sdppoint}. If $(x_+^\text{QCQP},x_-^\text{QCQP})$ is weakly $\epsilon$-efficient for the LP of robustness level $r^\text{QCQP}$, 
	then $(x_+^\text{QCQP},x_-^\text{QCQP},r^\text{QCQP})^T$ is weakly $\epsilon$-efficient for the QCQP \eqref{eq:myqcqp}.
\end{theorem}

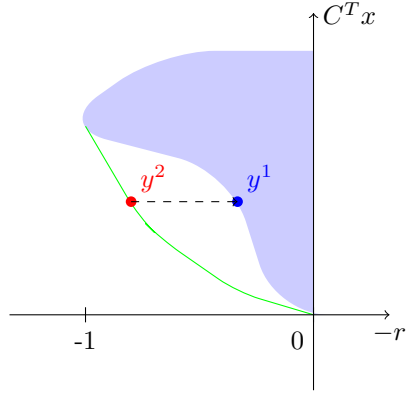
\begin{figure}
	\centering
	\begin{tikzpicture}[scale=1]

		\fill[blue!20,rounded corners=7mm] (3,3.5) -- (1,3.5) -- (-0.4,2.5) -- (1.9,1.9) -- (2.5,0) -- (3,0);		
		
		\draw[green,rounded corners=3mm] (0,2.5) --(0.7, 1.3) -- (1,1) -- (2, 0.3) -- (3,0);
		
		
		\fill[red] (0.6,1.5) circle (2pt) node[above right] {$y^2$};
		\fill[blue] (2,1.5) circle (2pt) node[above right] {$y^1$};
		
		\draw[->, dashed] (0.6,1.5) -- (2,1.5);
		
		\draw[->] (-1,0) -- (4,0) node[below] {$-r$};
		\draw[->] (3,-1) -- (3,4) node[right] {$C^Tx$};
		\draw[black] (3,0.1) -- (3,-0.1) node[below left] {0};
		\draw[black] (0,0.1) -- (0,-0.1) node[below] {-1};
	\end{tikzpicture}
	\caption{Visualization of the reconstruction method: The blue area indicates the image of the feasible set of the QCQP. Note that all components of $C^Tx$ have been combined into a single coordinate axis for this visualization. The green curve represents the Pareto-front of the SDP-relaxation with a Pareto-optimal point (red) on it. The blue point indicates its projection along the dashed line onto the QCQP and is, under the assumptions of Theorem \ref{thmoptstr}, again optimal.}
\end{figure}

Theorems \ref{thmoptstr} and \ref{thmoptw} offer a connection between optimal points of the QCQP \eqref{eq:myqcqp}, optimal points of the LP of robustness level $r^\text{QCQP}$ and feasible points of the SDP relaxation. Clearly, the condition on efficiency for the LP is a strong one in both theorems. The general idea of the projection onto the QCQP is that optimality for the SDP should roughly translate to optimality for the QCQP. 

The condition on efficiency for the LP serves mostly as a verification tool in a post-processing step, to ensure that an actually optimal point has been chosen. The method itself does not require to check this condition in order to reconstruct points that are feasible for the QCQP from the SDP relaxation (which in itself is remarkable). Further, the idea that optimality should translate from the SDP to the QCQP is not entirely unjustified, as the following chain of theorems shows:

\begin{theorem}[{\cite[Corollary 5.11]{schroeder}}]\label{coroptlb}
	Let $Z^\text{SDP}$ as in \eqref{eq:sdppoint} be a weakly efficient solution for the multi-criteria SDP. If $r^\text{SDP}=0$, then 
	\begin{align}
		\begin{pmatrix}
			x_+^\text{SDP}\\
			x_-^\text{SDP}\\
			r^\text{SDP}
		\end{pmatrix}
	\end{align}
	is weakly efficient for the QCQP \eqref{eq:myqcqp} with the same objective function value as the multi-criteria SDP.
\end{theorem}

A similar result holds for the other extreme case:

\begin{theorem}[{\cite[Corollary 5.13]{schroeder}}]\label{coroptub}
	Let $Z^\text{SDP}$ as in \eqref{eq:sdppoint} be a weakly efficient solution for the SDP. If $r^\text{SDP}=1$, then 
	\begin{align}
		\begin{pmatrix}
			x_+^\text{SDP}\\
			x_-^\text{SDP}\\
			r^\text{SDP}
		\end{pmatrix}
	\end{align}
	is weakly efficient for the QCQP \eqref{eq:myqcqp}.
\end{theorem}

In this case, a converse theorem is true as well:

\begin{theorem}[{\cite[Theorem 5.14]{schroeder}}]\label{thmoptub}
	Let $(x_+^\text{QCQP},x_-^\text{QCQP},r^\text{QCQP})^T$ be an efficient solution of the QCQP \eqref{eq:myqcqp} with $r^\text{QCQP}=1$. Then 
	\begin{align}
		\begin{pmatrix}1\\x_+^\text{QCQP}\\x_-^\text{QCQP}\\r^\text{QCQP}\end{pmatrix}
		\begin{pmatrix}1\\x_+^\text{QCQP}\\x_-^\text{QCQP}\\r^\text{QCQP}\end{pmatrix}^T
	\end{align}
	is efficient for the SDP.
\end{theorem}

\section{A New Navigation Method for the Pareto-Front}\label{paretonavmeth}
Based on the accumulated results in the previous sections, we propose a new method for navigation of the Pareto-front of \eqref{eq:myqcqp}. It consists of the following steps:

\begin{enumerate}
	\item Calculate the Pareto-front of the SDP relaxation.
	\item Navigate the convex Pareto-front of the SDP via convex navigation tools to a point $Z^\text{SDP}$.
	\item Project $Z^\text{SDP}$ onto the feasible set of \eqref{eq:myqcqp} via \eqref{eq:maxr}.
\end{enumerate}

Since step 3 can be done in real time, the decision-maker gets direct feedback for his currently navigated point. Once a suitable point $(x_+,x_-,r)^T$ has been found, we can run an automated check on its optimality with Theorems \ref{thmoptstr} or \ref{thmoptw}. 

If the point is already optimal, we are done. Otherwise, we can re-optimize within the slice of the selected robustness level $r$ to find a solution that is at least as good as the one selected by the decision-maker. Since $r$ remains fixed for this last re-optimization step, this is once again an LP that can be solved efficiently.

This method has the potential of being further enhanced by making use of Theorem \ref{coroptub} and Theorem \ref{thmoptub}. Since for $r=1$ the Pareto-fronts of the relaxation and the original QCQP coincide, a decision-maker could, for example, first explore the worst-case front, by fixing $r$ to $1$. All Pareto-points that are found in this way can then already be added to the Pareto-front of the SDP. This saves computation time when the inverse robust optimization via the SDP relaxation is started.

\section{Comparison with Established Methods}\label{compmeth}
In order to illustrate our method with a practical example, in this section we consider a large-scale liver case with many variables and constraints that serves as an example of practically relevant magnitude. We demonstrate the computability of solutions, compare the number of iterations, the runtime and the solution quality to those of the $\epsilon$-constraint method. Finally, we discuss observed advantages and disadvantages of our method.

We use the open source software for radiation treatment planning of intensity-modulated photon, proton, and carbon ion therapy MatRad \footnote{https://e0404.github.io/matRad/, see \cite{wieser}}. We work with its liver case, which consists of multiple structures and two objective functions: dose maximization on the PTV and dose minimization on the skin tissue. 

We exported the dose-influence matrix $D$ from MatRad and imposed the bounds from Table \ref{tab:dosebounds} on each structure. We replaced the quadratic objectives functions (squared overdose/deviation) with the linear mean dose objective function in the nominal case. Additionally, we set $\ell=0$ and $u=50$.

\begin{table}
	\begin{tabular}{|c|c|c|}
		\hline
		Structure & Lower Bound (\si{\gray}) & Upper Bound (\si{\gray})\\
		\hline
		PTV/CTV/GTV		 & 40 & 50\\
		Spinal Chord     & 0  & 20\\
		Heart			 & 0  & 40\\
		Heart (mean)     & 0  & 10\\
		Liver-CTV        & 0  & 50\\
		Others           & 0  & 45\\ 
		\hline 
	\end{tabular}
	\caption{Dose bounds for each structure.}
	\label{tab:dosebounds}
\end{table}

We converted the problem into the form $Ax\leq b$ with only one-sided inequality constraints, i.~e.,
\begin{align}
	A_c=\begin{pmatrix}
		D\\
		-D
	\end{pmatrix}
\end{align}	
and
\begin{align}
	b_c=\begin{pmatrix}
		b^\text{ub}\\
		-b^\text{lb}
	\end{pmatrix}.
\end{align}
Further we set up the problem to handle errors within the dose-influence matrix of up to $2\%$, i.~e.,
\begin{align}
	A_\delta=0.02\cdot A_c.
\end{align}
The right-hand side generally does not contain uncertainty, since we assume that the physician is confident about the prescribed lower and upper doses and so we set 
\begin{align}
	b_\delta=0.
\end{align}
Finally, for the mean objective functions we extracted the dose-influence sub-matrix $D^s$ for structure $s$ from $D$ and set 
\begin{align}
	g^s_j=\frac{1}{\vert I(s)\vert}\sum_{i=1}^{\vert I(s)\vert}D^s_{ij} \text{ for } j=1,2,\dots,n. 
\end{align}
We accumulate these objective functions to the objective matrix
\begin{align}
	G=\begin{pmatrix}
		(-g^\text{PTV})^T\\
		(g^\text{Skin})^T
	\end{pmatrix}.
\end{align}

The $-1$ in front of $g^\text{PTV}$ indicates that we want to maximize the dose on this structure. With these choices of parameters, we obtain a nominal LP with $923$ variables and (after clustering) $30813$ constraints in addition to the bound constraints on $x$. 

A similar approach with interval uncertainty for the IMRT problem has recently been proposed in \cite{moreno}. But there, in contrast to what we are doing, the authors follow the classical method of iteratively incrementing the robustness level $r$. 

After application of Theorem \ref{IntMatLinear} to our problem and introduction of the inverse robustness variable $r$ as well as the auxiliary variable $x_0=1$, we get the QCQP with $1848$ variables. Before we reformulate this problem into the SDP-relaxation, we make the following observation: From the bound $0\leq x=x_+-x_-$ we have $x_-\leq x_+$, so that whenever $x_->0$ we know that the pair $(x_+-x_-,0)$ is also feasible with the same or a better objective function value. Hence, we can set $x_-=0$ and remove these variables from the QCQP, reducing the problem size to $925$ variables. 

This has a significant impact, as the relaxation now only has a matrix of size $925\times 925$ instead of $1848\times 1848$, which amounts to roughly a factor of $4$ for the number of variables and consequently results in a significant reduction of memory usage and computation time. 

Now we formulate our LP constraints and our bounds, as well as the additionally introduced bounds on the diagonal of $Z$ and $Z\geq0$ as SDP constraints, i.~e., as linear constraints on the variable $Z$. Since adding the constraint $Z_{ij}\geq0$ for all $i$ and $j$ would lead once again to an extremely large problem, we only add those constraints that immediately interact with our data, i.~e., we only impose the constraints $Z_{i,1}\geq0$ and $Z_{925,j}\geq0$. Note that due to the symmetry of $Z$, this also implies $Z_{1,j}\geq0$ and $Z_{i,925}\geq0$. While the other entries of $Z$ surely are important too, we found them not to be as significant for our method, since they are mostly multiplied by the zeros of the constraint matrices $M_i$. They only play a role to make $Z$ positive semidefinite. 

We applied the sandwiching algorithm \cite{lammelsandw}, which is an algorithm for Pareto-front approximation, to this problem with an accuracy of $\delta=0.04$, which was achieved after calculating a total of $19$ Pareto-points. In contrast to this, we solved the QCQP with incremental increases in $r$ of $0.04$ and solved the resulting $26$ LPs with the same sandwiching algorithm to a quality of $\delta=0.04$. These sandwiching runs accumulated to a total of $163$ Pareto-points. 

We see that with this method roughly eight times as many Pareto-points had to be computed. This observation reinforces our conjecture that using the sandwiching algorithm on the SDP instead of iteratively solving the QCQP leads to much fewer computations of Pareto-points. 

On the other hand, Table \ref{tab:imrtruntimes} shows the computation times for both approaches and we observe much higher runtimes for our approach in this example. This can be attributed to the immense size of the SDP problem. Fine-tuning on the clustering method should allow to further reduce the runtimes to achieve practically useful ones. 

\begin{figure}
	\includegraphics[width=\textwidth]{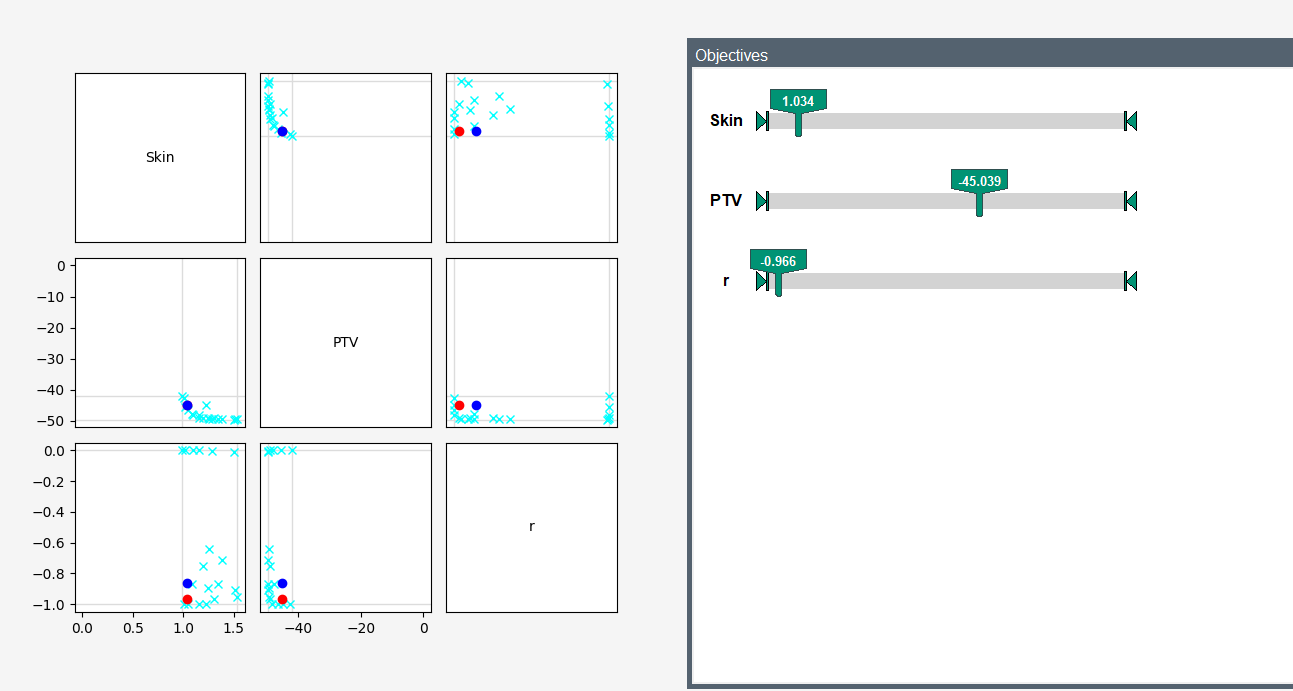}
	\caption{Pareto-front of the relaxed IMRT problem, visualized in Fraunhofer ITWMs Pareto navigation tool. The matrix on the left showcases the projection of all Pareto-points onto the 2D-plane with only the indicated functions. The lightblue crosses represent the Pareto-points. On the right, the values of the currently navigated point (red) of the SDP are shown. Further, its projection onto the QCQP (blue) is displayed. Its coordinates are $(1.034, -45.06, -0.862)$. A particular Pareto-point of the QCQP is given by $(1.0339, -45.0665, -0.84)$, making the point $\epsilon$-optimal for $\epsilon =0.022$.}
	\label{fig:navigation}
\end{figure}

Further, we confirmed manually in the navigator that the projection onto the QCQP finds $\epsilon$-optimal Pareto-points (with $\epsilon\approx0.022$). Figure \ref{fig:navigation} shows this in more detail.

\begin{table}
	\centering
	\begin{tabular}{|c|c|c|}
		\hline
		Method&Computed Pareto-points&Runtime\\\hline
		SDP-relaxation & 19 & 46h44m\\
		Iterative robustness levels & 163 & 2h31m\\
		\hline
	\end{tabular}
	\caption{Comparison of our SDP-relaxation and the method of iteratively incremented robustness levels $r$ for the IMRT problem.}
	\label{tab:imrtruntimes}
\end{table}

\begin{figure}
	\begin{subfigure}{\textwidth}
		\includegraphics[width=.5\textwidth]{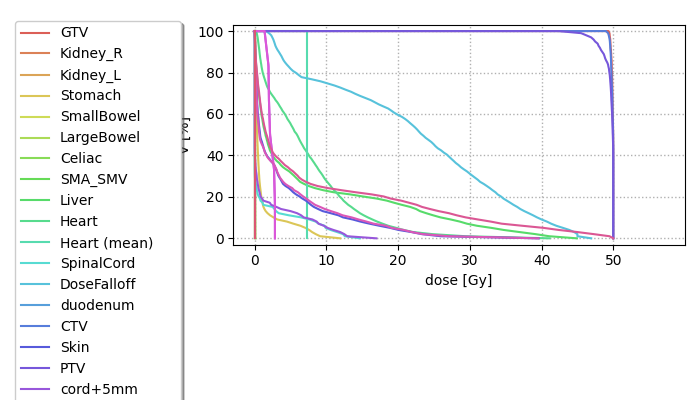}
		\includegraphics[width=.5\textwidth]{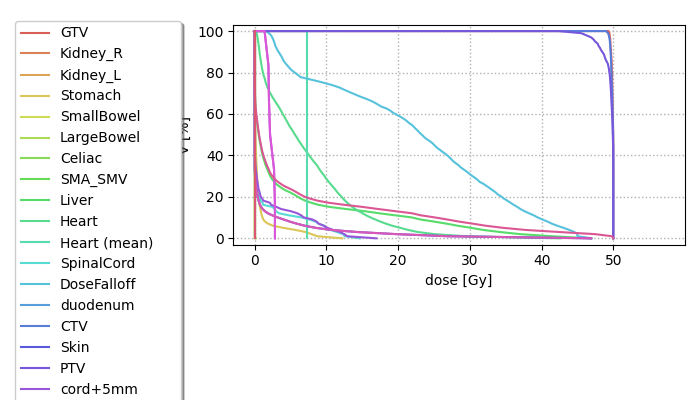}
		\caption{Objective: $(1.381, -49.544, -0.077)^T$}
		\label{fig:mydvhsr=0.08}
	\end{subfigure}
	\begin{subfigure}{\textwidth}
		\includegraphics[width=.5\textwidth]{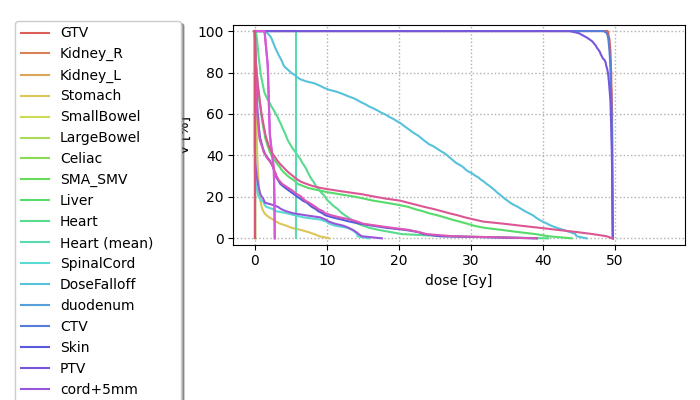}
		\includegraphics[width=.5\textwidth]{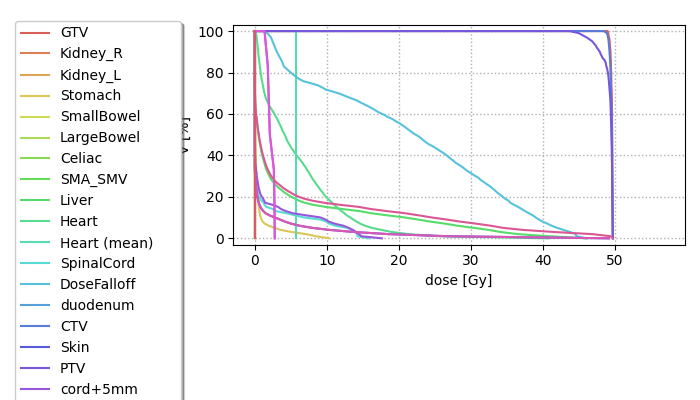}
		\caption{Objective: $(1.304, -49.106, -0.411)^T$}
		\label{fig:mydvhsr=0.4}
	\end{subfigure}
	\begin{subfigure}{\textwidth}
		\includegraphics[width=.5\textwidth]{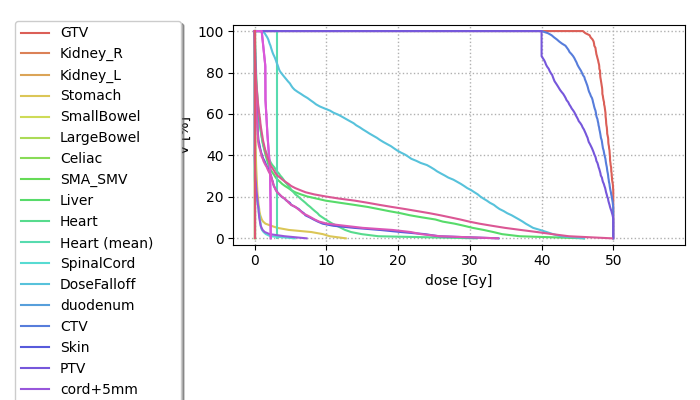}
		\includegraphics[width=.5\textwidth]{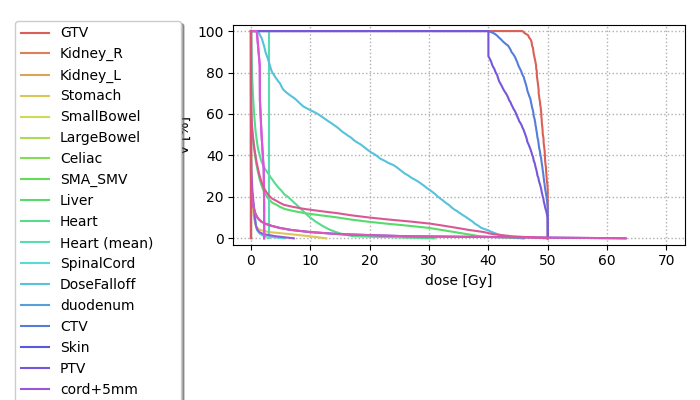}
		\caption{Objective: $(1.015, -45.488, -0.703)^T$}
		\label{fig:mydvhsr=0.7}
	\end{subfigure}
	\begin{subfigure}{\textwidth}
		\includegraphics[width=.5\textwidth]{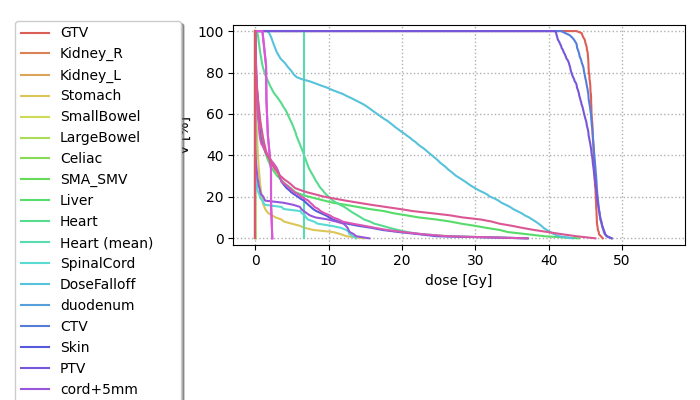}
		\includegraphics[width=.5\textwidth]{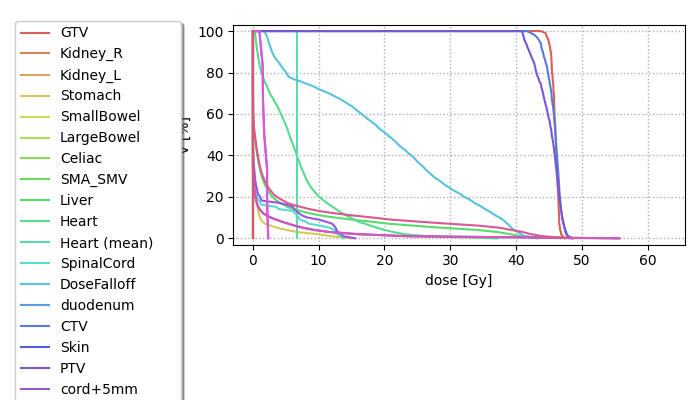}
		\caption{Objective: $(1.226, -44.973, -1.000)^T$}
		\label{fig:mydvhsr=1}
	\end{subfigure}
	
	\caption{DVHs at different Pareto-points of the QCQP. On the left are the DVHs for the clustered problems, on the right for the unclustered problem with the same intensity vector $x$. From Figure \ref{fig:mydvhsr=0.08} to \ref{fig:mydvhsr=0.4} we see a slight shift away from the upper bound of $50\si{\gray}$ to make some ``breathing room'' for the uncertainties to unfold. Figures \ref{fig:mydvhsr=0.4} and \ref{fig:mydvhsr=0.7} are not really comparable, because they represent very different Pareto-points, Figure \ref{fig:mydvhsr=0.4} putting much more emphasis on high tumor doses. From Figure \ref{fig:mydvhsr=0.7} to \ref{fig:mydvhsr=1}, note how with increasing levels of robustness the dose on the target becomes more homogenous (the curves representing the tumor doses are less spread in \ref{fig:mydvhsr=1} compared to \ref{fig:mydvhsr=0.7}).}
	\label{fig:mydvhs}
\end{figure}

We further evaluated our navigated points in terms of DVHs. Figure \ref{fig:mydvhs} showcases a selection of DVHs at Pareto-points of different levels of robustness. 

\begin{table}[t]
	\centering
	\begin{tabular}{|c|c|}
		\hline
		Pareto-point of the SDP & Projection onto the QCQP\\
		\hline
		$(1.381, -49.544, -0.712)^T$&$(1.381, -49.544, -0.0766)^T$\\
		$(1.304, -49.106, -0.968)^T$&$(1.304, -49.106, -0.411)^T$\\
		$(1.085, -47.685, -0.871)^T$&$(1.085, -47.685, -0.471)^T$\\
		$(1.226, -44.973, -1.0)^T$&$(1.226, -44.973, -1.0)^T$\\
		\hline
	\end{tabular}
	\caption{Selected Pareto-points of the SDP and their projection onto the QCQP.}
	\label{tab:paretoproj}
\end{table}

Finally, table \ref{tab:paretoproj} shows a selection of Pareto-points of the SDP and their projection onto the QCQP. We observe quite a big ``loss'' in terms of $r$, bringing many of the conceivably robust Pareto-points to much worse robustness levels. The interpretation here is that the Pareto-front is particularly steep for large values of $r$. 

This may be due to the immense dimension $p=925\cdot30813=28502025$ of our uncertainty set: Recall that for the utility function we chose $\vartheta(W)=r$. While this ensures that our problem remains linear, the function $\vartheta(W)=r^p$ would have closer resembled the volume of the uncertainty set. 

One could argue that our choice greatly distorts the perception of the number of scenarios: Even if $r=0.99$, we have $r^p<0.01$, i.~e. less than $1\%$ of the volume of the full uncertainty set is considered. This is a strictly numerical problem though and the interpretation that each individual parameter could assume up to $99\%$ of its maximum deviation will certainly still be relevant. Nevertheless, this could be one explanation for the large fall-off in the projection. 

On the other hand it highlights the importance of Theorems \ref{coroptub} and \ref{thmoptub}: By calculating the fully robust front (i.~e., $r=1$) separately, we can always navigate towards robust solutions and be sure to ultimately achieve the robustness that we desire.

\section{Conclusions}\label{conclusion}

Let us summarize the key results of this paper: Beginning  with a linear radiotherapy problem, we introduced uncertainty in the dose-influence matrix. We enhanced the classical worst-case approach with methods from interval arithmetic and the novel concept of inverse robustness, yielding a framework that is free of semi-infinite constraints and results in a (non-convex) QCQP. 

We relaxed this problem to obtain an SDP formulation, which can be efficiently solved with interior-point solvers. Clustering played a crucial role in reducing the problem size to make the problem computationally tractable. Next, starting from a solution of the SDP, we presented a reconstruction method to recover a solution to the QCQP. This method is backed by multiple theorems that ensure (weak) efficiency under certain conditions and theorems that relate the Pareto-front of the SDP and the QCQP to each other. 

Finally, we showed the applicability of this method for IMRT problems, showcasing a significant reduction in the number of calculated Pareto-points, while retaining the same approximation quality. While the runtime for very large problems is significantly higher than with the method of iteratively incrementing the robustness, in small and moderate dimensions there may be optimization problems, where our method is faster, too. 

Notably, both methods show a way to find solutions of QCQPs (which are generally NP-hard) in polynomial time. This is only possible due to the structure of the particular QCQPs. Moreover, the example showed how inverse robustness can be integrated into practical problems to obtain relevant solutions that are more robust than solutions of a nominal optimization. The modelling and the combination with the Pareto-navigation tool allows for an interactive exploration of the trade-offs between robustness and other objective functions. This allows decision-makers who are not mathematicians to get an intuitive understanding of the properties of a navigated point, allowing for quicker and more confident decisions.

\medskip{}
\medskip{}
\textbf{Statements and Declarations}

\textbf{Competing Interests}: The authors have no competing interests to declare that are relevant to the content of this article.

All authors contributed to the study conception and design. All authors read and approved the final manuscript.

Data sets generated during the current study are available from the corresponding author on reasonable request.

\textbf{Acknowledgements}: The work of Y.C. is supported by U.S. National Institutes of Health Grant Number R01CA266467 and by the Cooperation Program in Cancer Research of the German Cancer Research Center (DKFZ) and Israel's Ministry of Innovation, Science and Technology (MOST).

\nocite{*}
\bibliography{bibliography}

\end{document}